\documentclass[11pt]{article}
\usepackage[T1]{fontenc}
\usepackage{lmodern,amsmath,amssymb,amsthm,booktabs,needspace}
\usepackage[letterpaper,hmargin=1.1in,vmargin=1.0in]{geometry}
\usepackage[hidelinks]{hyperref}
\newtheorem{theorem}{Theorem}[section]
\newtheorem{proposition}[theorem]{Proposition}
\newtheorem{conjecture}{Conjecture}
\DeclareMathOperator{\Gal}{Gal}
\DeclareMathOperator{\IMG}{IMG}
\DeclareMathOperator{\St}{St}
\DeclareMathOperator{\Rist}{Rist}
\newcommand{\Q}{\mathbb Q}
\newcommand{\Z}{\mathbb Z}
\newcommand{\Pp}{\mathbb P^1}
\title{Sections and Constant Fields of\\Postcritically Finite Rational Maps}
\author{Wayne Peng\\Department of Mathematics, National Central University}
\hypersetup{pdftitle={Sections and Constant Fields of Postcritically Finite Rational Maps},pdfauthor={Wayne Peng}}
\date{}
\begin{document}
\maketitle
\begin{abstract}
We relate the constant fields of generic preimage towers to two invariants of geometric monodromy. Bounded section indices characterize branching of the closed action and imply a finite constant extension. Unbounded finite-level abelianization orders force an infinite cyclotomic subextension. For a quadratic rational map with hyperbolic orbifold, we compute the section indices and abelianizations at every sufficiently large level and obtain an infinite real cyclotomic subextension. We also give a branching criterion for quadratic maps with neither critical point periodic and formulate three conjectures about the discrete action, the constants, and the case of bounded abelianizations with unbounded section indices.
\end{abstract}

\section{The geometric action and the constants}

Let $f\in K(z)$ be a postcritically finite (PCF) rational map of degree $d\geq2$ over a number field $K$, and let $P=P_f=\bigcup_{i\geq1}f^{\circ i}(\operatorname{Crit}f)$. For transcendental $t$, write $E_n$ for the splitting field of $f^{\circ n}(z)=t$, $C_n=E_n\cap\overline K$, and $C_\infty=\bigcup_n C_n$. Over $\mathbb C$, monodromy on iterated preimages of a base point in $U=\Pp(\mathbb C)\setminus P$ gives $\Gamma=\IMG(f)$. Write $G$ for its tree closure and $G_n$ for its level-$n$ image. Comparison with geometric monodromy gives
\[
1\longrightarrow G_n\longrightarrow\Gal(E_n/K(t))
 \longrightarrow\Gal(C_n/K)\longrightarrow1.
\]
For a finite coefficient extension $F/K$, the splitting field becomes $FE_n$ and the constants become $FC_n$, by regularity over $C_n$. Geometric coordinate conjugacy preserves the tree action.

Known exact calculations often start with a recognizable recursion or a special critical portrait. This makes it hard to tell which patterns persist beyond those families. We use section indices and finite abelianizations to compare the examples. Earlier work includes the descriptions of quadratic polynomials by Bartholdi--Nekrashevych and Pink \cite{barth-nek,pink}, the constant fields of unicritical polynomials in Adams--Hyde \cite{adams-hyde}, and arithmetic and geometric models for quadratic rational portraits in Ejder and Ejder et al.\ \cite{ejder-periodic,ejder}. The general branch-to-constants implication is Fari\~na-Asategui's theorem \cite[Theorem 2]{farina}.

\subsection{Section indices}

Use the usual connecting-path labels for the tree \cite{nek-book}. The action is self-replicating: it is level transitive, and the section at a child maps its stabilizer onto $\Gamma$. Indeed, the covers are connected, and the section map comes from the surjection of fundamental groups induced by filling the extra punctures in $f^{-1}(U)\subseteq U$. The same holds for $G$. Write $\St_J(n)$ for the pointwise level stabilizer of a tree group $J$ and $\Rist_J(n)$ for the product of its rigid stabilizers at level $n$. A level-transitive action is branch when each $\Rist_J(n)$ has finite index.

Let $H=\psi(\St_G(1))\leq G^d$, with $\psi$ the section map, and let $H_n$ be its image in $G_n^d$. Its index
\begin{equation}\label{eq:index}
 I_n=[G_n^d:H_n]=\frac{|G_1||G_n|^d}{|G_{n+1}|}
\end{equation}
measures the constraints on the child sections. The order formula comes from the kernel of $G_{n+1}\to G_1$.

\begin{proposition}\label{prop:index}
The $I_n$ are nondecreasing, and $\sup_n I_n=[G^d:H]$ as an extended integer. They are bounded exactly when $G$ is branch. Bounded $I_n$ implies $[C_\infty:K]<\infty$.
\end{proposition}
\begin{proof}
The finite coset maps are surjective, and closedness of $H$ makes any finite set of distinct cosets distinguishable at some common level. A bounded index therefore makes $H$ open. Then $H$ contains $\St_G(m)^d$ for some $m$, and inserting these sections gives parents in $\St_G(m)$ itself. Iterating gives independent copies at every level, proving branching. Conversely, the section image of $\Rist_G(1)$ is a product of coordinate subgroups of finite index in the respective projections of $H$. Those projections have finite index in $G$ by self-replication, since the full root stabilizer has finite index in each child stabilizer. The constant-field claim is \cite[Theorem 2]{farina}.
\end{proof}

The discrete section image $H_\Gamma=\psi(\St_\Gamma(1))$ is dense in $H$. The same finite-index argument works discretely: if $H_\Gamma$ has finite index, its coordinate intersections have finite index, and successive section maps show that the section image of each full level stabilizer has finite index in the full section product. Its coordinate intersections then give the rigid stabilizers. The converse uses the projections as above. In degree two, both section projections are onto. Goursat's lemma identifies their image as a fiber product over a common quotient by the subgroup that can be inserted in either child. The two insertion kernels coincide by normality and a first-level swap.

\paragraph{Finite certificates.}
For a finitely generated action, one can certify $[\Gamma^d:H_\Gamma]\leq q$ using $q$ proposed right cosets $Uc_i$, where $U$ is generated by verified section tuples. Include $c_1=1$ and check $c_i s\in Uc_{j(i,s)}$ for a symmetric generating set of $\Gamma^d$. The union is then the whole group. Conversely Schreier's lemma supplies such certificates for finite index. With finite-state generators the word equalities can be checked on reachable sections.

There is also an exact finite recurrence when a finite quotient $\Lambda\to Q$ of a self-similar action has kernel $N$ satisfying $N^d\leq\psi(N\cap\St_\Lambda(1))$. In a closed action take $N$ open. Let $J\leq Q^{d+1}$ record the quotient label of each root-fixing parent and its sections, and let $T$ be its child projection. Starting at $U_0=Q$, set $U_{n+1}$ equal to the parent projection of $J\cap(Q\times U_n^d)$. Then $U_n$ is the image of $\St_\Lambda(n)$ and $I_n=[Q^d:TU_n^d]$. Indeed, each child label in $U_n$ can be corrected by $N$ to a level-$n$ stabilizer element, with no change to the parent label. Equality of consecutive $U_n$ makes the recurrence constant. An initial plateau in $I_n$ alone gives no such certificate. For related algorithms for regularly branched actions see Bartholdi--Noce \cite{barth-noce}.

\subsection{Inertia and abelianizations}
Let $B_n=G_n^{\mathrm{ab}}$, $b_n=|B_n|$, and $e_n=\exp B_n$. The peripheral loops about $P$ generate $G_n$ and can be taken with one loop omitted, so $b_n\leq e_n^{|P|-1}$. The order of a peripheral loop in $\Gamma$ is its ramification weight $\nu(p)$, the least common multiple of all local degrees of iterates over $p$. In particular $b_n$ is bounded when at most one weight is infinite.

\begin{theorem}\label{thm:cycl}
Put $L=K(P)$ and $M_n=K(\mu_{e_n})$. Then $M_n\subseteq LC_n$ and
\[
 [C_n:K]\geq [L(\mu_{e_n}):L]\geq \frac{\varphi(e_n)}{[L:\Q]}.
\]
If $b_n$ is unbounded, then $C_\infty\cap K(\mu_\infty)$ is infinite over $K$.
\end{theorem}
\begin{proof}
Over $LC_n$, the punctures are rational and the arithmetic group of the level cover is $G_n$. The branch-cycle lemma \cite[Lemma 6.4]{adams-hyde} sends each inertia class to its cyclotomic power. Arithmetic lifts act by inner conjugation on $G_n$, so on $B_n$ the cyclotomic power fixes its generators. Thus it is the identity modulo $e_n$. This proves the containment and the degree bound by taking degrees over $L$. Also $[M_n:M_n\cap C_n]\leq[L:K]$, so $[M_n\cap C_n:K]\geq\varphi(e_n)/[L:\Q]$. If $b_n$ is unbounded, so are $e_n$ and $\varphi(e_n)$.
\end{proof}

\section{Three regimes and a hyperbolic example}

The two criteria leave the mixed regime, where $I_n$ is unbounded and $b_n$ bounded. The table gives representatives over $\Q$. Write $\Q(\mu_{2^\infty})^+$ for the union of the real $2$-power cyclotomic fields. The Chebyshev formulas hold from level two, and those for Latt\`es and $f_*$ from level three.
\begin{center}
\small
\begin{tabular}{@{}llll@{}}
\toprule
$f$ & $I_n$ & $b_n$ & $C_\infty$ or a subfield\\
\midrule
$z^2$ & $2^n$ & $2^n$ & $\Q(\mu_{2^\infty})$\\
$z^2-1$ (Basilica) & $2^{\lfloor n/2\rfloor}$ & $2^n$ & $\Q(\mu_{2^\infty})$\\
$z^2-2$ (Chebyshev) & $2^{n+1}$ & $4$ & $\Q(\mu_{2^\infty})^+$\\
$2/(z-1)^2$ (Latt\`es) & $2^{n+2}$ & $8$ & contains $\Q(\mu_{2^\infty})$\\
$z^2-16/(27z)$ (Hanoi) & $2$ & $2$ & $\Q(\sqrt3)$\\
$f_*$ in Theorem~\ref{thm:mixed} & $2^{1+\lfloor n/2\rfloor}$ & $8$ & contains $\Q(\mu_{2^\infty})^+$\\
\bottomrule
\end{tabular}
\end{center}

The power map follows from $z^{2^n}=t$. For Basilica, Pink gives the two cyclic factors of the abelianization and the diagonal section quotient of order $2^{\lfloor n/2\rfloor}$ \cite[Lemmas 2.2.5 and 2.2.9, Theorem 2.8.2]{pink}. Adams--Hyde compute its constants \cite[Theorem 1.6(2)]{adams-hyde}. For Chebyshev, $t=u^m+u^{-m}$ with $m=2^n$ gives roots $\zeta_m^j u+\zeta_m^{-j}u^{-1}$. Rotation and inversion give the geometric dihedral action, and simultaneous inversion of $u,\zeta_m$ gives the real constant field \cite[Theorem 1.6(3), case (D)]{adams-hyde}. Here $I_1=1$ and $b_1=2$. The Latt\`es row uses \cite[Lemma 3.3 and Theorems 1.6--1.9]{ejder}. In this case $C_\infty/\Q$ is nonabelian.

\paragraph{The Hanoi calculation.}
Grigorchuk--\v{S}uni\'c identify the geometric tree action with the standard three-peg action and give discrete commutator subgroup $D$ with $[\Gamma:D]=8$ and $[D:D^3]=12$ under section insertion \cite[Example 8 and Theorem 5.1]{gs}. Since $G_1=S_3$, the discrete section index is $8^3/(8\cdot12/6)=32$ in that model. Radi identifies $\overline D$ with the root-sign kernel in $G$ \cite[Proposition 5.20]{radi}, giving $b_n=2$. With $a=(a,1,1)(12)$, $b=(1,b,1)(02)$, $c=(1,1,c)(01)$ and composition from right to left, every root-fixing section tuple has even sign sum. The section image contains $\overline D^3$, and $abac=(a,a,bc)$ and its conjugates give all even sign patterns. Thus $I_n=2$.

At a vertex $v$ the preimage cubic is $z^3-vz-16/27$, with discriminant $\Delta(v)=4v^3-256/27$. For its roots $w$, $\Delta(w)=4vw-64/9$, so $\prod_w\Delta(w)=(1024/27)\Delta(v)$. The four Vandermondes at $v$ and its children multiply to $\pm(32/9)\sqrt3\,\Delta(v)$. This gives $\sqrt3\in E_2$. Order each fiber recursively to make every sign positive. In these labels the arithmetic action has sign discrepancy $s(v)+\sum_{f(w)=v}s(w)$ equal to the character of $\sqrt3$ at every vertex with two more levels below it. The group with zero discrepancies has orders $M_1=6$, $M_{n+1}=3M_n^3$, since the child actions determine the root sign. These equal $|G_n|$ by \eqref{eq:index}. Hence $C_1=\Q$ and $C_n=\Q(\sqrt3)$ for $n\geq2$.

\begin{theorem}\label{thm:mixed}
Let $f_*(z)=(-3z^2+4z-1)/(4z-1)$. Its orbifold is hyperbolic. It has $I_1=I_2=1$ and, for $n\geq3$, writing $m_n=2^{\lfloor n/2\rfloor}$,
\[
 I_n=2m_n,\qquad G_n^{\mathrm{ab}}\simeq(\Z/2\Z)^3,\qquad
 \Q(\zeta_{m_n}+\zeta_{m_n}^{-1})\subseteq C_n.
\]
\end{theorem}
\begin{proof}
The critical orbits are $0\leftrightarrow1$ and $1/2\mapsto1/4\mapsto\infty\mapsto\infty$, with $\infty$ unramified. Thus the weights are $(\infty,\infty,2,2)$. Appendix~\ref{app:calc} gives the marked recursion and the all-level calculation. Its surviving common section quotient is dihedral of order $2m_n$, and the appendix supplies a quotient of $G_n$ onto this dihedral group in which the rational meridian at $0$ is a generating rotation. By the branch-cycle lemma, a lift of an automorphism fixing $C_n$ takes the meridian to a conjugate of its cyclotomic power. The lift acts by inner conjugation on $G_n$, so the two powers are conjugate in $G_n$ itself. In the dihedral quotient a generating rotation is conjugate precisely to its two signed powers. The cyclotomic exponent is consequently $\pm1$ modulo $m_n$, giving the stated containment.
\end{proof}

The example gives maps of degree $2^k$ in the mixed regime by taking $k$th iterates. More generally, the level groups and fields for $f^{\circ k}$ are $G_{kn}$ and $E_{kn}$. Iterating \eqref{eq:index} gives
\[
 I^{(k)}_n=\prod_{j=0}^{k-1}(I_{kn+j}/I_j)^{d^{k-1-j}},\qquad I_0=1 .
\]
Boundedness of each sequence and the total constant field are unchanged. The orbifold weights agree as well, since any backward path can be extended at its start to a length divisible by $k$.

\paragraph{Scope of the data.}
The versioned ADD cohort contains 97 displayed coordinate records \cite{add}. For its 59 hyperbolic records the recorded classifications are 21 trivial, 24 finite nontrivial, and 14 infinite. The other 38 are parabolic and all are recorded infinite. These are arithmetic classifications over the respective coefficient fields of the displayed models. Eleven records were revised in the study, and 86 retain baseline classifications. The ancillary files \texttt{add-constant-tower-records.csv} and \texttt{add-kh15-action-bounds.csv} give these records and 15 separate action-level bounds from the November 22 release. The latter are upper bounds for the specified stored actions. The two extracts quantify the sample behind the examples.

\section{A branching criterion in degree two}

\begin{theorem}\label{thm:quad}
Suppose neither critical point of a quadratic PCF rational map is periodic. Let $j$ be the deck involution and $C=P\cap j(P)$. If $C\ne P$, both the discrete and closed actions are branch. If $C=P$, the map is Latt\`es and both actions are virtually abelian and nonbranch.
\end{theorem}
\begin{proof}
The common quotient of the two discrete sections is a quotient of a Euclidean orbifold group on $C$. The proof uses an entering meridian whose lift has support in one child, whereas its image in an infinite affine quotient would be a rotation.

Indeed, the section maps $q_0,q_1$ from $\pi_1(\Pp(\mathbb C)\setminus f^{-1}P)$ are onto $\Gamma$. They fill outside $P$ and $j(P)$ respectively, with the second map transported by $j$. If $K_i=\ker q_i$, the common quotient is the quotient by $K_0 K_1$. Thus points outside $C$ disappear, and each remaining meridian has order dividing its weight. If $r$ critical points belong to $P$, then $r\leq1$. The graph on $P$ has $2-r$ sources, the critical values whose critical preimages are outside $P$, and hence $2-r$ vertices with two distinct predecessors. Their pairs of predecessors, together with the $r$ marked critical points, give $C$. For $r=0$ all weights are two. For $r=1$ the unique marked critical point has weight two and the points after it have weight four. The two predecessors of the cycle-entry vertex occur after that critical point, since a critical value cannot be periodic under the hypothesis. The signatures on $C$ are therefore $(2,2,2,2)$ and $(2,4,4)$.

The corresponding groups are $\Z^2\rtimes R$ with $R$ generated by inversion or a quarter-turn. In any infinite quotient the kernel lies in the translation lattice, since a nontrivial linear part in the kernel has commutators spanning a full-rank lattice. Saturating the kernel gives an infinite affine quotient with torsion-free translations in which all nontrivial rotations act without fixed vectors and all cone generators retain nontrivial linear part. A first-child rigid motion must map to a translation there. Otherwise its commutators with the translations from the image of the root stabilizer generate a full-rank translation subgroup in its rigid image. A root swap transports this subgroup into the other rigid image. The rigid images commute, forcing the original linear part to fix a full-rank lattice, a contradiction.

Each cycle in $P$ contains a point of $C$: the periodic predecessor and the incoming predecessor at the first entry into the cycle are paired by $j$. The entry cannot be a critical value. Thus if $C\ne P$ some $p\in P\setminus C$ maps to $q\in C$. The two preimages of $q$ are $p$ and the unmarked $j(p)$, so the meridian at $q$ has one trivial section and is rigid. It still represents a cone generator, which contradicts the infinite affine quotient. The common section quotient is therefore finite, giving discrete and closed branching.

The overlap contributes two to $\sum_{p\in P}(1-1/\nu(p))$, so the full orbifold is hyperbolic for $C\ne P$ and Euclidean for $C=P$. In the latter case the pullback weights $\nu(f(x))/\deg_x f$ dominate $\nu(x)$ and have the same Euler characteristic by Riemann--Hurwitz. They therefore equal $\nu(x)$, giving the flat orbifold covering of a Latt\`es map \cite[Theorem 4.1]{milnor}. Both groups are virtually abelian since the Euclidean orbifold group surjects onto $\Gamma$. An infinite level-transitive virtually abelian tree group is nonbranch: at a sufficiently deep level, a finite-index normal abelian subgroup moves every vertex. Each rigid stabilizer then intersects it trivially by disjoint supports and so is finite.
\end{proof}

\section{Three conjectures}

\begin{conjecture}
If the tree closure of the IMG of a PCF rational map is branch, then its discrete IMG is branch.
\end{conjecture}
When $H$ is open, $H\cap\Gamma^d$ has index $[G^d:H]$ in $\Gamma^d$. The conjecture asks whether $H_\Gamma$ has finite index in this intersection. Hanoi gives indices $32$ and $2$ for $H_\Gamma$ and $H\cap\Gamma^d$ in the displayed model. Theorem~\ref{thm:quad} proves the conjecture in its quadratic range.

Which critical portraits or section quotients characterize the mixed regime? Theorem~\ref{thm:quad} gives exactly the Latt\`es maps in its range, while Theorem~\ref{thm:mixed} provides a hyperbolic family under iteration. The following arithmetic conjecture concerns what the known families share.

\begin{conjecture}\label{conj:abel}
If $I_n$ is unbounded and $b_n$ bounded, then for some finite $F/K$ the extension $FC_\infty/F$ contains an infinite abelian extension of $F$.
\end{conjecture}

\begin{conjecture}\label{conj:fin}
For a PCF rational map over a number field, $[C_\infty:K]<\infty$ if and only if $I_n$ is bounded.
\end{conjecture}

Proposition~\ref{prop:index} gives one direction of Conjecture~\ref{conj:fin}. The reverse would follow from Conjecture~\ref{conj:abel}: for unbounded $I_n$, unbounded $b_n$ gives infinite constants by Theorem~\ref{thm:cycl}, and bounded $b_n$ gives infinite constants after finite base change by the conjecture. Compare \cite[Conjecture 3]{farina} for hyperbolic-orbifold polynomials with nontrivial constants.

\begin{proposition}\label{prop:exception}
Every signed Chebyshev or Latt\`es map is in the mixed regime and satisfies the arithmetic assertion of Conjecture~\ref{conj:abel}.
\end{proposition}
\begin{proof}
For signed Chebyshev maps, the level-$n$ iterate is, up to sign, the Chebyshev map of degree $d^n$. The dihedral calculation above applies, and the constants contain the real cyclotomic fields \cite[Theorem 2.6]{hamblen-jones}. These conclusions persist after finite extension to a field of coordinate conjugacy.

For Latt\`es maps, take an elliptic cyclic quotient $\pi:X\to\Pp$ of degree $q$ and an affine lift $\psi=\alpha+Q$ of $f$ \cite[Theorem 3.1]{milnor}. These data descend to a finite $F/K$: the cyclic cover has algebraic branch values, so it descends by its Kummer equation, and the graph of the lift is a component of the descended fiber product $\pi(B')=f(\pi(B))$. The deck group fixes a chosen origin. Over algebraically closed constants the roots are $\pi(B+T)$, with $t=\pi(\psi^n B)$ and $T\in A_n=\ker\alpha^n$. The $d^n$ images are distinct since a nontrivial deck rotation identifying two lifts would make $(h-1)B$ constant for a nonzero endomorphism $h-1$. Thus $d^n\leq |G_n|\leq qd^n$, proving unbounded $I_n$. The orbifold has finite weights, giving bounded $b_n$.

Put $D_n=F(A_n)$. Lifting $\pi(B)$ and the images of $B+T_i$ for two generators of $A_n$ costs degree at most $q^3$ over $FE_n$. Subtracting $B$ recovers the torsion generators. In particular $[D_n:D_n\cap FC_n]\leq q^3$. This is the lifting argument of Hamblen--Jones \cite[Theorem 2.7]{hamblen-jones}, applied to an arbitrary isogeny kernel. If $\alpha=[m]$, take the cyclotomic subfield $F(\mu_{|m|^n})$ of $D_n$ from the Weil pairing. Its intersection with $FC_n$ has unbounded degree. Otherwise, on each rational Tate module the nonscalar $\alpha$ has commutative centralizer, so the primary torsion fields over $F$ are abelian. Since $A_n\subseteq X[d^n]$ by the dual isogeny, $D_n/F$ is abelian. The $D_n$ form a tower of unbounded degree by finiteness of elliptic torsion over a number field. Their intersections with $FC_n$ therefore give the claimed infinite abelian extension.
\end{proof}

\paragraph{Acknowledgment.}
OpenAI Codex assisted with drafting, mathematical checks, and editing.

\appendix
\section{The calculation for the hyperbolic example}\label{app:calc}

\paragraph{The marked action.}
Order the real intervals cut by $0,1/4,1,\infty$ from left to right and let $r_1,\ldots,r_4$ be their reflections in the quotient by complex conjugation. In Nekrashevych's real-reflection construction \cite{nek-reflections}, the sheet over the upper half-plane inside $|z-1/4|=1/4$ lies above the real line and the exterior sheet below it. This follows from
\[
 \operatorname{Im} f_*(z)=\frac{3\,\operatorname{Im}z}{4}\,
 \frac{1/16-|z-1/4|^2}{|z-1/4|^2}.
\]
The inverse discriminant is $4(t-1)(4t-1)$. For $t=-1,1/8,2$ the inverse pairs are $4/3\pm\sqrt{10}/3$, $7/12\pm\sqrt7/12$, and $(-2\pm\sqrt7)/3$. Reading their intervals in inside--outside order gives $r_1=(r_3,r_4)$, $r_2=(r_3,r_3)$, $r_4=(r_2,r_1)$. On $(1/4,1)$ the lifts meet along the circle, so $r_3=\sigma$ with trivial sections. A section word $w$ from sheet $i$ to sheet $j$ is corrected to $r_3^iwr_3^j$. The adjacent products $A=r_1r_2,B=r_2r_3,C=r_3r_4$ give the sphere recursion, in right-action order,
\begin{equation}\label{eq:marked}
 A=(1,C),\qquad B=\sigma,\qquad C=(AB,B)\sigma .
\end{equation}
In particular $A$ is a meridian at $0$ up to orientation.

\paragraph{The common quotient.}
Use the formal group $\mathcal F=\langle A,B,C\mid B^2=1\rangle$ and denote its level kernels by $K_n$. Its root stabilizer has generators $A,BAB,CB,BC$, with section pairs $(1,C),(C,1),(AB,B),(B,AB)$. The projections are onto. The first two pairs and the squares of the last two give isolated copies of $C,(AB)^2$, so the section image contains $N^2$ for their normal closure $N$. Modulo $N$ it is the graph interchanging the generating reflections $X=\eta(AB)$ and $Y=\eta(B)$ of $D_\infty=\mathcal F/N$, where $\eta$ is the quotient map and $R=\eta(A)=XY$. Write $\iota$ for this swap. With $M_n=\eta(K_n)$ the finite section relation is the fiber product over $D_\infty/(M_n\iota(M_n))$, with the swap on one side. Indeed two representatives can be adjusted by $K_n$ exactly when the discrepancy in the dihedral group lies in $M_n\iota(M_n)$. Thus $I_n$ is the order of this quotient. We will find $M_n=\langle R^{m_n}\rangle$ for $n\geq3$.

\paragraph{The level kernels.}
Use an auxiliary formal group $F_0=\langle a,b,z\mid a^2=b^2=1\rangle$ with recursion
\[
 a=(bazb,bz^{-1}ab)\sigma,\qquad b=(bab,b),\qquad
 z=(ba,bz^{-1}b)\sigma .
\]
For $g^T=T^{-1}gT$, substitution with $T=(bT,bz^{-1}aT)$ gives $A=(azb)^T$, $B=a^T$, $C=(az^{-1}a)^T$ with recursion \eqref{eq:marked}. Conversely $a,z,b$ are recovered from these three words. Thus the map from $\mathcal F$ to $F_0$ given by these words before conjugation is onto and carries level kernels onto level kernels. In what follows write $K_n^0$ for those of $F_0$. The homomorphism $\eta$ factors through $F_0$, taking $a,b,z$ to $Y,YXY,1$. These two reflections multiply to $R$.

In $H=\St_{F_0}(1)$ put $t=za,u=az,v=b,w=aba$, so $H=\langle t,u,v,w\mid v^2=w^2=1\rangle$. Let $L$ be the kernel of the map sending all four generators to $1\in C_2$. It is free on $(tv,uv,vt,vu,wv)$ by Schreier's lemma. In $H/[L,L]$ write a label as $(\mathbf n,e)\in\Z^5\rtimes C_2$, using $v$ for the odd representative. The involution swaps the first and third coordinates, swaps the second and fourth, and negates the fifth. It is convenient to use the unimodular coordinates
\[
 (x,\delta,r,p,q)=(n_1-n_3+n_5,\ n_1+n_2-n_3-n_4+2n_5,\ n_5,n_3,n_4).
\]
In these coordinates the involution is $J(x,\delta,r,p,q)=(-x,-\delta,-r,x-r+p,-x+\delta-r+q)$. The labels of $t,u,v,w$ have sign one and vectors $(1,1,0,0,0),(0,1,0,0,0),0,(1,2,1,0,0)$, respectively.

Here are the section pairs of a Schreier generating set of $K_2^0$ in $H$. Capitals indicate inverses, and the transversal for the two child parities is $1,t,v,w$.
\begin{center}
\small
\begin{tabular}{@{}rlll@{\qquad}rlll@{}}
\toprule
& Parent & Section $0$ & Section $1$ & & Parent & Section $0$ & Section $1$\\
\midrule
1 & $uT$ & $vtwv$ & $vUwTuv$ & 7 & $vuw$ & $vuvUv$ & $UwTuv$\\
2 & $tt$ & $vwTvUv$ & $1$ & 8 & $vwT$ & $vtv$ & $v$\\
3 & $tu$ & $vwTv$ & $vU$ & 9 & $wtv$ & $vuvUwTv$ & $v$\\
4 & $tvw$ & $vwTuvUv$ & $vUtwTuv$ & 10 & $wuv$ & $vuvUv$ & $vUv$\\
5 & $twv$ & $vTv$ & $v$ & 11 & $wvT$ & $vuvUtwv$ & $vUtwTuv$\\
6 & $vtw$ & $vUv$ & $v$ & & & &\\
\bottomrule
\end{tabular}
\end{center}

Both section projections are onto $H$. For instance, the product of rows $1,2,6^{-1}$ gives $v$ on the left, and row 5 gives $v$ on the right. Rows $5,6,3$ then supply the other left generators. On the right, rows $3,7,4$ supply $u,wT,twT$, giving $t,w$. With $Y=(1,1)$ and $YXY=(0,1)$ in dihedral rotation/sign coordinates, the map on $H$ is $(\delta,e)$. The rows obey $\delta_1=\delta_0+e$ with equal signs. The left kernel contains the normal closure of $u^{-1}t,t^2,w t v t$: use rows $6\cdot 5^{-1}$, $8\cdot 5^{-1}$, and row 2 modulo the first two relations. The presentation after imposing these relations is dihedral. Consequently the section relation is exactly the set of pairs with matching signs and $\delta_1=\delta_0+e$. Surjectivity on the right lets one lift any right entry and correct the left entry by that kernel.

On these pairs of labels the parent vector is
\[
 \bigl(-\delta_0-e,\ 2(x_0+x_1)+e,\ 2x_0+e,\
 x_0+r_0-p_0-q_0+e,\ -\delta_0-x_1+r_1-p_1-q_1-e\bigr)
\]
with sign $e$. To verify the formula, multiply labels by $(h,e)(k,f)=(h+J^e k,e+f\bmod2)$. Substitution shows the parent rule is a homomorphism on the stated pairs, and applying it to the table gives the labels of its parent words. This proves the rule for the whole section image.

For the even labels of $K_n^0$, the three moduli of $(x,\delta,r)$ start at $(1,1,1)$ for $n=1$. The rule gives the exact update $(a_0,b_0,c_0)\mapsto(b_0,2a_0,2a_0)$, since $x_0,x_1$ vary independently and the last two coordinates remain unrestricted. At level two every label has $\delta\equiv e\pmod2$. Two such labels with matching signs cannot realize $\delta_1=\delta_0+e$ when $e=1$, so all labels from level three on are even. The rotation modulus is then $2^{\lfloor n/2\rfloor}$, as required. The $z$-exponent of an even label is $\delta+2p+2q-2r$. The $a$-parity on $H$ is the parity of this exponent, and the $b$-parity is $e$ minus that parity. Thus $K_n^0$ maps to exactly $(0,0,2\Z)$ in $F_0^{\mathrm{ab}}=C_2(a)\times C_2(b)\times\Z(z)$ for $n\geq3$, proving $G_n^{\mathrm{ab}}\simeq C_2^3$. Finally $M_1=D_\infty$, and $M_2=\langle R^2,Y\rangle$: even rotation exponents at level two are exactly the even integers, and odd labels occur with odd rotation exponent by the parent rule. The swap gives $\iota(M_2)=\langle R^2,X\rangle$, so $M_2\iota(M_2)=D_\infty$ and $I_1=I_2=1$.

\end{document}